\documentclass{article}

\usepackage{amsmath}
\usepackage{amsfonts}
\usepackage{amssymb}
\usepackage{amscd}
\usepackage{amsbsy}
\usepackage{amsthm}
\usepackage{graphicx}
\usepackage{psfrag}
\usepackage{color}

\usepackage{cleveref}
\usepackage{mathrsfs}

\usepackage[twoside]{geometry}
\newfont{\myfnt}{cmssi10 scaled 1440}
\numberwithin{equation}{section}
\catcode`@=11
\def\ps@nk{\def\@oddhead{\vbox{\hbox to \hsize{\pic \footnotesize \it \shorttitle
				\hfill \rm \thepage} \vspace{1mm} \vspace*{-2mm}}}
	\def\@evenhead{\vbox{\hbox to \hsize{\pic \footnotesize \rm \thepage \hfill \it \shortauthor}
			\vspace{1mm} \vspace*{-2mm}}}
	\def\@oddfoot{} \def\@evenfoot{}}
\def\ps@first{\def\@oddhead{}
	\def\@evenhead{}
	\def\@oddfoot{} \def\@evenfoot{}}
\def\ps@total{\def\@oddhead{\vbox{\hbox to \hsize{\footnotesize \rm \hfill TOTAL\ \ CONTENTS
				\hfill \thepage} \vspace{1mm} \hrule \vspace*{-2mm}}}
	\def\@evenhead{\vbox{\hbox to \hsize{\footnotesize \rm \thepage \hfill CHIN.\ \ ANN.\ \ MATH.
				\hfill} \vspace{1mm} \hrule \vspace*{-2mm}}}
	\def\@oddfoot{} \def\@evenfoot{}}
\newtheoremstyle{thmstyle}
{6pt}
{6pt}
{\it}
{}
{\bf}
{}
{.5em}
{}
\newtheoremstyle{remstyle}
{6pt}
{6pt}
{\rm}
{}
{\bf}
{}
{.5em}
{}

\def\Section#1{\Sec{\large #1} \setcounter{equation}{0} \vskip -6mm \indent}
\def\Sec{\@Startsection{section}{1}{\z@}
	{-3.5ex \@plus -1ex \@minus -.2ex}%
	{2.3ex \@plus.2ex}%
	{\normalfont\large\bfseries\boldmath}}
\def\@Startsection#1#2#3#4#5#6{%
	\if@noskipsec \leavevmode \fi
	\par
	\@tempskipa #4\relax
	\@afterindenttrue
	\ifdim \@tempskipa <\z@
	\@tempskipa -\@tempskipa \@afterindentfalse
	\fi
	\if@nobreak
	\everypar{}%
	\else
	\addpenalty\@secpenalty\addvspace\@tempskipa
	\fi
	\@ifstar
	{\@ssect{#3}{#4}{#5}{#6}}%
	{\@dblarg{\@Sect{#1}{#2}{#3}{#4}{#5}{#6}}}}
\def\@Sect#1#2#3#4#5#6[#7]#8{%
	\ifnum #2>\c@secnumdepth
	\let\@svsec\@empty
	\else
	\refstepcounter{#1}%
	\protected@edef\@svsec{\@seccntformat{#1}\relax}%
	\fi
	\@tempskipa #5\relax
	\ifdim \@tempskipa>\z@
	\begingroup
	#6{%
		\@hangfrom{\hskip #3\relax\@svsec \hskip -2.5mm}%
		\interlinepenalty \@M #8\@@par}
	\endgroup
	\csname #1mark\endcsname{#7}%
	\addcontentsline{toc}{#1}{%
		\ifnum #2>\c@secnumdepth \else
		\protect\numberline{\csname the#1\endcsname}%
		\fi
		#7}%
	\else
	\def\@svsechd{%
		#6{\hskip #3\relax
			\@svsec #8}%
		\csname #1mark\endcsname{#7}%
		\addcontentsline{toc}{#1}{%
			\ifnum #2>\c@secnumdepth \else
			\protect\numberline{\csname the#1\endcsname}%
			\fi
			#7}}%
	\fi
	\@xsect{#5}}
\renewenvironment{abstract}{%
	\small
	\quotation
	\noindent {\bfseries \abstractname } }%
{\if@twocolumn\else\endquotation\fi}

\def\Subsec{\@StartSubsection{subsection}{2}{\z@}%
	{-3.25ex\@plus -1ex \@minus -.2ex}%
	{1.5ex \@plus .2ex}%
	{\normalfont\normalsize\bfseries\boldmath}}
\def\@StartSubsection#1#2#3#4#5#6{%
	\if@noskipsec \leavevmode \fi
	\par
	\@tempskipa #4\relax
	\@afterindenttrue
	\ifdim \@tempskipa <\z@
	\@tempskipa -\@tempskipa \@afterindentfalse
	\fi
	\if@nobreak
	\everypar{}%
	\else
	\addpenalty\@secpenalty\addvspace\@tempskipa
	\fi
	\@ifstar
	{\@ssect{#3}{#4}{#5}{#6}}%
	{\@dblarg{\@SubSect{#1}{#2}{#3}{#4}{#5}{#6}}}}
\def\@SubSect#1#2#3#4#5#6[#7]#8{%
	\ifnum #2>\c@secnumdepth
	\let\@svsec\@empty
	\else
	\refstepcounter{#1}%
	\protected@edef\@svsec{\@seccntformat{#1}\relax}%
	\fi
	\@tempskipa #5\relax
	\ifdim \@tempskipa>\z@
	\begingroup
	#6{%
		\@hangfrom{\hskip #3\relax\@svsec\hskip -1.5mm}%
		\interlinepenalty \@M #8\@@par}
	\endgroup
	\csname #1mark\endcsname{#7}%
	\addcontentsline{toc}{#1}{%
		\ifnum #2>\c@secnumdepth \else
		\protect\numberline{\csname the#1\endcsname}%
		\fi
		#7}%
	\else
	\def\@svsechd{%
		#6{\hskip #3\relax
			\@svsec #8}%
		\csname #1mark\endcsname{#7}%
		\addcontentsline{toc}{#1}{%
			\ifnum #2>\c@secnumdepth \else
			\protect\numberline{\csname the#1\endcsname}%
			\fi
			#7}}%
	\fi
	\@xsect{#5}}
\def\list#1#2{\ifnum \@listdepth >5\relax \@toodeep \else \global
	\advance \@listdepth\@ne \fi \rightmargin \z@ \listparindent\z@
	\itemindent\z@ \csname @list\romannumeral\the\@listdepth\endcsname
	\def\@itemlabel{#1}\let\makelabel\@mklab \@nmbrlistfalse #2\relax
	\@trivlist \parskip 0pt \parindent\listparindent \advance \linewidth
	-\rightmargin \advance\linewidth -\leftmargin \advance\@totalleftmargin
	\leftmargin \parshape \@ne \@totalleftmargin \linewidth \ignorespaces}
\renewcommand{\@makecaption}[2]{\begin{center}#1. #2\end{center}}\catcode`@=12 

\theoremstyle{thmstyle}
\newtheorem{thm}{\indent Theorem}[section]
\newtheorem{lem}{\indent Lemma}[section]
\newtheorem{prop}{\indent Proposition}[section]

\theoremstyle{remstyle}
\newtheorem{rem}{\indent \bf Remark}[section]
\newtheorem{ex}{\indent \bf Example}[section]

\newsavebox{\mygraphic}

\def\pic{\begin{picture}(0,0) \put(-210,-1250){\usebox{\mygraphic}} \end{picture}}
\newfont{\HUGEbf}{cmbx10 scaled 3500}
\definecolor{gray}{rgb}{0.9,0.9,0.9}
\def\thebibliography#1{\section*{\bf \large References}
	\list{[\arabic{enumi}]} {\settowidth \labelwidth{[#1]} \leftmargin
		\labelwidth \advance \leftmargin \labelsep \usecounter{enumi}}
	\def\newblock{\hskip .11em plus .33em minus .07em} \footnotesize \sloppy \clubpenalty
	4000 \widowpenalty 4000 \sfcode`\.=1000 \relax}

\newcommand{\To}{\rightarrow}
\newcommand{\as}{{\rm d}\mathbb{P}\times{\rm d} t-a.e.}

\newcommand{\ps}{\mathbb{P}-a.s.}

\newcommand{\F}{\mathcal{F}}
\newcommand{\E}{\mathbb{E}}

\newcommand{\s}{\mathcal{S}}

\newcommand{\hcal}{\mathcal{H}}

\newcommand{\T}{[0,T]}

\newcommand{\R}{{\mathbb R}}

\newcommand {\Dis}{\displaystyle}

\def\firstpage{1}

\def\shorttitle{Multi-dimensional Diagonally Quadratic BSDEs} 
\def\shortauthor{{\it G. Yang}} 

\title{\Large \bf \boldmath\ \\ Multi-dimensional Backward Stochastic Differential Equations of Diagonally Quadratic Generators with a Special Structure$^{\ast}$}

\author{\large  Guang YANG$^1$} 

\date{}

\begin{document}
	\maketitle
	\thispagestyle{first}
	\renewcommand{\thefootnote}{\fnsymbol{footnote}}
	\footnotetext{\hspace*{-5mm} \begin{tabular}{@{}r@{}p{13.4cm}@{}}
$^1$ & Shanghai Center for Mathematical Sciences, Fudan University, Shanghai 200433, China.\\
&{E-mail: gyang19@fudan.edu.cn} \\
$^{\ast}$ & Research partially supported by National Natural Science Foundation of China (Grants No. 11631004 and No. 12031009).
\end{tabular}}
\renewcommand{\thefootnote}{\arabic{footnote}}
	
\begin{abstract}
The present paper is devoted to the well-posedness of a type of multi-dimensional backward stochastic differential equations (BSDEs) with a diagonally quadratic generator. We give a new priori estimate, and prove that the BSDE admits a unique solution on a given interval when the generator has a sufficiently small growth of the off-diagonal elements (i.e., for each $i$, the $i$-th component of the generator has a small growth of the $j$-th row $z^j$ of the variable $z$ for each $j \neq i$). Finally, we give a solvability result when the diagonally quadratic generator is triangular.
\vskip 4.5mm
		
		\noindent \begin{tabular}{@{}l@{ }p{10.1cm}} {\bf Keywords } &
			Multi-dimensional BSDE, diagonally quadratic generator, BMO martingale
		\end{tabular}
		
		\noindent {\bf 2000 MR Subject Classification } 
		60H10
	\end{abstract}
	
	\baselineskip 14pt
	
	\setlength{\parindent}{1.5em}
	
	\setcounter{section}{0}
	
\Section{Introduction} \label{section1}

Bismut~\cite{Bismut1973JMAA} first introduced Backward stochastic differential equations (BSDEs in short):
\begin{equation}\label{eq:1.1}
Y_t=\xi+\int_t^T f(s,Y_s,Z_s){\rm d}s-\int_t^T Z_s {\rm d}W_s, \ \ t\in\T,
\end{equation}
where $(W_t)_{t\in\T}$ is a $d$-dimensional standard Brownian motion defined on some complete probability space $(\Omega, \F, \mathbb{P})$, and $(\F_t)_{t\in\T}$ is the augmented natural filtration generated by the standard Brownian motion $W$. The terminal value $\xi$ is an $\F_T$-measurable $n$-dimensional random vector, the generator function $f(\omega, t, y, z):\Omega\times\T\times\R^n\times\R^{n\times d}\To \R^n$
is $(\F_t)$-progressively measurable for each pair $(y,z)$, and the solution $(Y_t,Z_t)_{t\in\T}$ is a pair of $(\F_t)$-progressively measurable processes with values in $\R^n\times\R^{n\times d}$ which almost surely verifies BSDE \eqref{eq:1.1}.
In 1990, Pardoux and Peng~\cite{PardouxPeng1990SCL} established the existence and uniqueness result for BSDEs with an $L^2$-terminal value and a generator satisfying a uniformly Lipschitz continuous condition. When the generators have a quadratic growth in the state variable $z$, the situation is more complicated. In the one-dimensional case, Kobylanski~\cite{Kobylanski2000AP} established the first existence and uniqueness result for quadratic BSDEs with bounded terminal values, Tevzadze~\cite{Tevzadze2008SPA} gives a fixed-point argument, Briand and Elie~\cite{BriandElie2013SPA} provide a constructive approach to quadratic BSDEs with and without delay. Briand and Hu~\cite{BriandHu2006PTRF,BriandHu2008PTRF},
Delbaen et al.~\cite{DelbaenHuRichou2011AIHPPS,DelbaenHuRichou2015DCD}, Barrieu and El Karoui~\cite{BarrieuElKaroui2013AoP} and Fan et al.~\cite{FanHuTang2019ArXiv} considered the unbounded terminal value case.

For multidimensional quadratic BSDEs, when the terminal value is small enough in the supremum norm, Tevzadze~\cite{Tevzadze2008SPA} proved a general existence and uniqueness result for multi-dimensional quadratic BSDEs. Frei and Dos Reis~\cite{FreiDosReis2011MFE} provide a counterexample which show that multidimensional quadratic BSDEs with a bounded terminal value may fail to have a global solution. Frei~\cite{Frei2014SPA} introduced the notion of split solution and studied the existence of solution by considering a special kind of terminal value. Cheridito and Nam~\cite{cheridito2015Stochastics} and Xing and {\v{Z}}itkovi{\'c}~\cite{XingZitkovic2018AoP} obtained the solvability for multidimensional quadratic BSDEs in the Markovian setting. Jamneshan et al.~\cite{JamneshanKupperLuo2017ECP} provided solutions for multidimensional quadratic BSDEs with separated generators.
Cheridito and Nam~\cite{cheridito2015Stochastics}, Hu and Tang~\cite{HuTang2016SPA} and Luo~\cite{Luo2019ArXiv} obtained local solvability of sysetems of BSDEs with subquadratic, diagonally quadratic and triangularly quadratic generators respectively, which under additional assumptions on the generator can be extended to global solutions. When the terminal value is unbounded, Jamneshan et al.~\cite{JamneshanKupperLuo2017ECP} provided solutions when the terminal value is small in the BMO-sense, Fan et al.~\cite{FanHuTang2020ArXiv} obtained global solutions when the generator is convex or concave.

As a continuation of Hu and Tang~\cite{HuTang2016SPA} and Fan et al.~\cite{FanHuTang2020ArXiv}, we are devoted to the solvability of multidimensional diagonally quadratic BSDEs when the generator has a small growth of the off-diagonal elements. The local solution is constructed directly by \cite[Theorem 2.2 ]{HuTang2016SPA}. Together with the new priori estimate we build and a special kind of `intermediate value' property of the $\s^{\infty}$-norm of the local solution, we are able to stitch local solutions to get the global solution. In contrast to \cite[Theorem 2.3]{HuTang2016SPA} and \cite[Theorem 2.4]{FanHuTang2020ArXiv}, we allow the generator to have a small growth of the off-diagonal elements. In contrast to \cite[Theorem 2.5]{FanHuTang2020ArXiv} and \cite{Luo2019ArXiv}, we do not assume that the generator is strictly quadratic. Finally, assuming that for each $i=1,\cdots,n$, the $i$th component $f^i$ of the generator $f$ is diagonally quadratic, depends only on the first $i$ components of the state variable $y$ and the first $i$ rows of the state variable $z$, we prove existence and uniqueness of the global solution to the multidimensional diagonally quadratic BSDE with a bounded terminal value.

The rest of the paper is organized as follows. In Section 2, we prepare some notations and state the main results of this paper. In Section 3, we give an estimate and prove our main results. In Section 4, we prove a global solvability result for triangular and diagonally quadratic BSDEs.




\Section{Preliminaries and statement of main results}
\subsection{Notations}
\vspace{0.1cm}

Let $W=(W_t)_{t\geq 0}$ be a $d$-dimensional standard Brownian motion defined on a complete probability space $(\Omega, \F, \mathbb{P})$, and $(\F_t)_{t\geq 0}$ be the augmented natural filtration generated by $W$. Throughout this paper, we fix a $T \in (0, \infty)$. We endow $\Omega \times \T$ with the predictable $\sigma$-algebra $\mathcal{P}$ and $\R^n$ with its Borel $\sigma$-algebra $\mathcal{B}(\R^n)$. All the processes are assumed to be $(\F_t)_{t\in\T}$-progressively measurable, and
all equalities and inequalities between
random variables and processes are understood in the sense of $\ps$ and $\as$,  respectively. The Euclidean norm is always denoted by $|\cdot|$,  and $\|\cdot\|_{\infty}$ denotes the $L^{\infty}$-norm for one-dimensional or multidimensional random variable defined on the probability space $(\Omega, \F, \mathbb{P})$.

We define  the following four Banach spaces of stochastic processes.
By $\s^p(\R^n)$ for $p\geq 1$ , we denote the set of all $\R^n$-valued continuous adapted processes $(Y_t)_{t\in\T}$ such that
$$\|Y\|_{{\s}^p}:=\left(\E[\sup_{t\in\T} |Y_t|^p]\right)^{1/p}<+\infty.$$
By $\s^{\infty}(\R^n)$, we denote the set of all $\R^n$-valued continuous adapted processes $(Y_t)_{t\in\T}$ such that
$$\|Y\|_{{\s}^{\infty}}:=\left\|\sup_{t\in\T} |Y_t| \right\|_{\infty}<+\infty.$$
By $\hcal^p(\R^{n\times d})$ for $p\geq 1$, we denote the set of all $\R^{n\times d}$-valued $(\F_t)_{t\in\T}$-progressively measurable processes $(Z_t)_{t\in\T}$ such that
$$
\|Z\|_{\hcal^p}:=\left\{\E\left[\left(\int_0^T |Z_s|^2{\rm d}s\right)^{p/2}\right] \right\}^{1/p}<+\infty.
$$
By ${\rm BMO}(\R^{n\times d})$, we denote the set of all $Z\in \hcal^2(\R^{n\times d})$ such that
$$
\|Z\|_{\rm BMO}:=\sup_{\tau}\left\|\E_{\tau}\left[\int_{\tau}^T |Z_s|^2 {\rm d}s\right]\right\|_{\infty}^{1/2}<+\infty.
$$
Here and hereafter the supremum is taken over all $(\F_t)$-stopping times $\tau$ with values in $\T$, and $\E_{\tau}$ denotes the conditional expectation with respect to $\F_\tau$.

The spaces $\s^p_{[a,b]}(\R^n)$, $\s^{\infty}_{[a,b]}(\R^n)$, $\hcal^p_{[a,b]}(\R^{n\times d})$,  and ${\rm BMO}_{[a,b]}(\R^{n\times d})$ are identically defined for stochastic processes over the time interval $[a,b]$. We note that for $Z\in {\rm BMO}(\R^{n\times d})$, the process $\int_0^t Z_s{\rm d}B_s, t\in\T$,  is an $n$-dimensional BMO martingale.  For the theory of BMO martingales, we refer the reader to Kazamaki~\cite{Kazamaki1994book}.

For $i=1,\cdots, n$, denote by $z^i$, $y^i$  and $f^i$ the $i$th row of matrix $z\in\R^{n\times d}$,  the $i$th component of the vector $y\in \R^n$ and the generator $f$, respectively.

\subsection{Statement of the main results}

The main result of this paper concerns global solutions for bounded terminal value case. Consider the multi-dimensional BSDE \eqref{eq:1.1} of the following structured
quadratic generator:
\begin{equation}\label{eq:2.1}
f^{i}(t,y,z)=g^{i}(t,z^{i})+h^{i}(t,y,z), \ \ i=1,\cdots, n.
\end{equation}

We need the following assumptions.

\begin{enumerate}
\renewcommand{\theenumi}{(H\arabic{enumi})}
\renewcommand{\labelenumi}{\theenumi}
\item\label{A:H1} There exist two positive real constants $\gamma$ and $C$ and a real constant $\delta \in [0,1)$, such that for $i=1,\cdots,n$, $g^i:\Omega\times\T\times\R^{d}\To \R$ and $h^i:\Omega\times\T\times\R^n\times\R^{n\times d}\To \R$ satisfy the following inequalities:
    $$
    \begin{array}{ll}
    &\Dis |g^i(t,z)|\leq \frac{\gamma}{2} |z|^2, \ \ \forall ~z \in \R^{d}; \vspace{0.1cm}\\
    &\Dis |g^i(t,z_1)-g^i(t,z_2)|\leq C(1+|z_1|+|z_2|)|z_1-z_2|, \ \ \forall ~z_1, z_2 \in \R^{d}; \vspace{0.1cm}\\
    &\Dis |h^i(t,0,0)|\leq C; \vspace{0.1cm}\\
    &\Dis |h^i(t,y_1,z_1)-h^i(t,y_2,z_2)|\leq C|y_1-y_2|+C(1+|z_1|^{\delta}+|z_2|^{\delta})|z_1-z_2|,\ \ \forall ~y_1,y_2 \in \R^n, \\
    &\Dis  z_1, z_2 \in \R^{n \times d}.
    \end{array}
    $$
\item\label{A:H2} There exist a three-dimensional non-negative deterministic vector function $(\alpha_{t},\beta_{t},\eta_{t})_{t\in\T}$ and a positive constant $r \in (0,1+\delta]$ such that for $i=1,\cdots,n$, $h^i$ satisfies:
    $$
    {\rm sgn}(y^i)h^i(t,y,z)\leq \alpha_{t}+\beta_{t}|y|+\eta_{t}|z|^{r},\ \ \forall ~y \in \R^n, z \in \R^{n \times d}.
    $$

\item\label{A:H3} There exist four non-negative constants $C_0$, $C_1$, $C_2$ and $C_3$ such that
$$\|\xi\|_{\infty}\leq C_0,\ \int_0^T \alpha_t{\rm d}t\leq C_1, \ \int_0^T \beta_t{\rm d}t\leq C_2 \ \ {\rm and}\ \int_0^T (\eta_t+\eta_{t}^{\frac{2}{1-\delta}}){\rm d}t \leq C_3.$$

\end{enumerate}

Our main result ensures existence and uniqueness for the diagonally quadratic BSDE\eqref{eq:1.1}.

\begin{thm}\label{thm:2.1}
There exists a constant $r_0>0$ (depending only on the vector of parameters $(n,\gamma,\delta,C_0,C_1,C_2,C_3)$) such that if \ref{A:H1}-\ref{A:H3} holds for $r \in (0,r_0)$, then BSDE \eqref{eq:1.1} has a unique solution $(Y,Z)\in \s^{\infty}(\R^n) \times {\rm BMO}(\R^{n\times d})$ on $[0,T]$.
\end{thm}
The proof is given in Section 3.
\begin{ex}\label{ex:2.1}
Assuming that $T=1$, then the following generator $f$ satisfies \ref{A:H1}-\ref{A:H3} with $(\alpha_{t},\beta_{t},\eta_{t})=(2,1,1)$ and $(\gamma,\delta,C_1,C_2,C_3)=(2,0.5,2,1,2)$ when $r \in (0,1.5]$:
$$f^i(t,y,z)=|z^i|^2+|y|+\sin(|z|^{\frac{3}{2}})+|z|^{r}{\bf 1}_{\{|z|>1\}}+|z|{\bf 1}_{\{|z|\leq 1 \}},\ \ i=1,\cdots, n.$$
\end{ex}

The second result of this paper concerns a special type of diagonally quadratic BSDEs as follows:
\begin{equation}\label{eq:2.2}
\Dis Y_{t}^{i}=\xi^{i}+\int_{t}^{T}k^{i}(s,Y_{s},Z_{s}){\rm d}s-\int_{t}^{T}Z_{s}^{i}{\rm d}W_s,\ 0 \leq t \leq T, \ 1 \leq i \leq n.
\end{equation}

For each $i=1,\cdots,n$, $H\in \R^{n\times d}$, $z\in \R^{1\times d}$, $Y\in \R^{n}$ and $y \in \R$, define by $H(z;i)$ the matrix in $\R^{n\times d}$ whose $i$th row is $z$ and whose $j$th row is $H^j$ for any $j\neq i$, define by $Y(y;i)$ the vector in $\R^{n}$ whose $i$th component is $y$ and whose $j$th component is $Y^j$ for any $j\neq i$. We make the following assumptions.
\begin{enumerate}
\renewcommand{\theenumi}{(A\arabic{enumi})}
\renewcommand{\labelenumi}{\theenumi}
\item\label{A:A1} There exist a constant $\alpha \in (-1,1)$ and a positive constant $K_1$ such that for $i=1,\cdots,n$, the function $k^i:\Omega\times\T\times\R^{n}\times\R^{n\times d}\To \R$ depends only on the first $i$ components of $y$ and the first $i$ rows of $z$, and:
    $$
    |k^i(t,y,z)| \leq K_{1}(1+\sum_{j=1}^{i}|y^j|+\sum_{j=1}^{i-1}|z^j|^{1+\alpha}+|z^i|^{2}),\ \ \forall ~y \in \R^{n},z \in \R^{n \times d}.
    $$
\item\label{A:A2} There exist a non-negative constant $\beta$ and a positive constant $K_2$ such that for $i=1,\cdots,n$ and each $(Y,Z,y_1,y_2,z_1,z_2)\in \R^{n}\times\R^{n\times d}\times \R \times \R \times \R^{1\times d} \times \R^{1\times d}$, the function $k^i$ satisfies:
    $$
    \Big|k^i\Big(t,Y(y_1;i),Z(z_1;i)\Big)-k^i\Big(t,Y(y_2;i),Z(z_2;i)\Big)\Big| \leq \beta\big|y_{1}-y_{2}\big|+K_2\big(1+|z_1|+|z_2|\big)\big|z_1-z_2\big|.
    $$

\item\label{A:A3} There exists a non-negative constant $K_3$ such that $\xi=(\xi^{1},\cdots,\xi^{n})^{*}$ satisfies
$$\|\xi\|_{\infty}\leq K_3.$$

\end{enumerate}

We have the following result.

\begin{thm}\label{thm:2.2}
Let Assumptions \ref{A:A1}-\ref{A:A3} be satisfied. Then BSDE \eqref{eq:2.2} has a unique solution $(Y,Z)\in \s^{\infty}(\R^n) \times {\rm BMO}(\R^{n\times d})$ on $[0,T]$.
\end{thm}
The proof is given in Section 4.

\begin{rem}\label{rmk2.1}
In Theorem~\ref{thm:2.2}, we do not require the assumptions \ref{A:H1}-\ref{A:H3}.
\end{rem}

\begin{ex}\label{ex:2.2}
The following generator $k$ satisfies \ref{A:A1}-\ref{A:A2} in Theorem~\ref{thm:2.2}:
$$
\begin{array}{lll}
& \Dis k^1(t,y,z)=1+y^1+\sin(y^1)+|z^1|^2;\\
& \Dis k^i(t,y,z)=1+\sum_{j=1}^{i}y^j+\sin(y^{i-1})y^i+\sum_{j=1}^{i-1}|z^j|^{1+\alpha}+\cos(|z^{i-1}|)|z^i|^{2},\ \ i=2,\cdots, n.
\end{array}
$$
\end{ex}

\Section{Diagonally quadratic BSDEs}
\vspace{0mm}

We first give an estimate. 

\begin{lem}\label{lem:3.1}
Let assumptions \ref{A:H1}-\ref{A:H3} hold, $(Y,Z) \in \s^{\infty}_{[t_0,T]}(\R^n) \times {\hcal^2}_{[t_0,T]}(\R^{n\times d})$ is a solution of BSDE \eqref{eq:1.1} on $[t_0,T]$, then there exist two positive constants $C_4,C_5$ (depending on the vector of parameters $(n,\gamma,\delta,C_0,C_1,C_2,C_3)$) such that

\begin{equation}\label{eq:3.1}
\|Y\|_{\s^{\infty}_{[t_0,T]}} \leq C_4+C_5\exp\big(\frac{r\gamma}{1-\delta}\|Y\|_{\s^{\infty}_{[t_0,T]}}\big).\vspace{0.2cm}
\end{equation}
\end{lem} 

{\bf Proof}\ \ Define
$$
u(x)=\frac{\exp(\gamma|x|)-\gamma|x|-1}{\gamma^2}, \ \ x \in \R.
$$
Then we have for $x \in \R$,
$$
u'(x)=\frac{\exp(\gamma|x|)-1}{\gamma}{\rm sgn}(x), \ \ u''(x)=\exp(\gamma|x|), \ \ u''(x)-\gamma |u'(x)| =1.
$$
Using It\^{o}'s formula to compute $u(Y_{t}^{i})$ and using the assumption \ref{A:H2}, we have
\begin{equation}\label{eq:3.2}
\begin{array}{lll}
\Dis u(Y_{t}^{i})& = &\Dis u(\xi^{i})+\int_{t}^{T}[u'(Y_{s}^{i})\big(g^{i}(s,Z_{s}^{i})+h^{i}(s,Y_{s},Z_{s})\big)-\frac{1}{2}u''(Y_{s}^{i})|Z_{s}^{i}|^{2}] {\rm d}s -\int_{t}^{T}u'(Y_{s}^{i})Z_{s}^{i}{\rm d}W_s\vspace{0.2cm}\\
& \leq & \Dis  u(\xi^{i})-\int_{t}^{T}u'(Y_{s}^{i})Z_{s}^{i}{\rm d}W_s \vspace{0.2cm}\\
& & \Dis +\int_{t}^{T}\Big[\frac{\exp(\gamma|Y_{s}^{i}|)-1}{\gamma}\big(\frac{\gamma}{2}|Z_{s}^{i}|^{2}+\alpha_{s}+\beta_{s}|Y_{s}|+\eta_{s}|Z_{s}|^{r}\big)-\frac{1}{2}\exp(\gamma|Y_{s}^{i}|)|Z_{s}^{i}|^{2}\Big] {\rm d}s \vspace{0.2cm}\\
& = & \Dis  u(\xi^{i})-\int_{t}^{T}u'(Y_{s}^{i})Z_{s}^{i}{\rm d}W_s \vspace{0.2cm}\\
& & \Dis +\int_{t}^{T}\Big[-\frac{1}{2}|Z_{s}^{i}|^{2}+\frac{\exp(\gamma|Y_{s}^{i}|)-1}{\gamma}\big(\alpha_{s}+\beta_{s}|Y_{s}|+\eta_{s}|Z_{s}|^{r}\big)\Big] {\rm d}s
\end{array}
\end{equation}
Using H\"{o}lder's inequality, we get
\begin{equation}\label{eq:3.3}
\eta_{s}|Z_{s}|^{r} = \varepsilon^{\frac{r}{2}}|Z_{s}|^{r} \cdot \varepsilon^{-\frac{r}{2}}\eta_{s} \leq \frac{r}{2}\varepsilon|Z_{s}|^{2}+ \frac{2-r}{2}(\varepsilon^{-\frac{r}{2}}\eta_{s})^{\frac{2}{2-r}}. 
\end{equation}
Taking 
$$
\varepsilon=\frac{\gamma}{nr}\exp(-\gamma\|Y\|_{\s^{\infty}_{[s,T]}}).
$$ 
We have
\begin{equation}\label{eq:3.4}
\Dis \eta_{s}|Z_{s}|^{r} \leq \frac{\gamma}{2n}\exp(-\gamma\|Y\|_{\s^{\infty}_{[s,T]}})|Z_{s}|^{2}+ \frac{2-r}{2}\eta_{s}^{\frac{2}{2-r}}\big(\frac{\gamma}{nr}\exp(-\gamma\|Y\|_{\s^{\infty}_{[s,T]}})\big)^{-\frac{r}{2-r}}.
\end{equation}
From $0< r \leq 1+\delta <2$, we have
\begin{equation*}
\frac{2-r}{2} \leq 1, \ \ 1 \leq \frac{2}{2-r} \leq \frac{2}{1-\delta}, \ \ \frac{r}{2-r} \leq\frac{1+\delta}{1-\delta}. 
\end{equation*}
Therefore
\begin{equation*}
\eta_{s}^{\frac{2}{2-r}} \leq \eta_{s}^{\frac{2}{1-\delta}}+\eta_{s}, \ \ (\frac{nr}{\gamma})^{\frac{r}{2-r}} \leq (\frac{2n}{\gamma})^{\frac{r}{2-r}} \leq (\frac{2n}{\gamma})^{\frac{1+\delta}{1-\delta}}+1.
\end{equation*}
From (\ref{eq:3.4}), we deduce that 
\begin{equation}\label{eq:3.5}
\Dis \eta_{s}|Z_{s}|^{r} \leq \Dis \frac{\gamma}{2n}\exp(-\gamma\|Y\|_{\s^{\infty}_{[s,T]}})|Z_{s}|^{2} + (\eta_{s}^{\frac{2}{1-\delta}}+\eta_{s})((\frac{2n}{\gamma})^{\frac{1+\delta}{1-\delta}}+1)\exp(\frac{r\gamma}{2-r}\|Y\|_{\s^{\infty}_{[s,T]}}).
\end{equation}
Let
$$
k_s=(\eta_{s}^{\frac{2}{1-\delta}}+\eta_{s})((\frac{2n}{\gamma})^{\frac{1+\delta}{1-\delta}}+1).
$$ 
From (\ref{eq:3.2}) and (\ref{eq:3.5}), we have 
\begin{equation}\label{eq:3.6}
\begin{array}{lll}
\Dis u(Y_{t}^{i}) & \leq & \Dis  u(\xi^{i})-\int_{t}^{T}u'(Y_{s}^{i})Z_{s}^{i}{\rm d}W_s+\int_{t}^{T}\big[-\frac{1}{2}|Z_{s}^{i}|^{2}+\frac{1}{2n}|Z_{s}|^{2}\big] {\rm d}s \vspace{0.2cm}\\
& & \Dis +\int_{t}^{T}\frac{\exp(\gamma|Y_{s}^{i}|)}{\gamma}\big(\alpha_{s}+\beta_{s}\|Y\|_{\s^{\infty}_{[s,T]}}+k_s\exp(\frac{r\gamma}{2-r}\|Y\|_{\s^{\infty}_{[s,T]}})\big) {\rm d}s.
\end{array}
\end{equation}
Hence it holds that
\begin{equation}\label{eq:3.7}
\begin{array}{lll}
\Dis \sum_{i=1}^{n}u(Y_{t}^{i}) & \leq & \Dis  \sum_{i=1}^{n}u(\xi^{i})-\int_{t}^{T} \sum_{i=1}^{n}u'(Y_{s}^{i})Z_{s}^{i}{\rm d}W_s \vspace{0.2cm}\\
& & \Dis +\frac{1}{\gamma}\int_{t}^{T}\big(\alpha_{s}+\beta_{s}\|Y\|_{\s^{\infty}_{[s,T]}}+k_s\exp(\frac{r\gamma}{2-r}\|Y\|_{\s^{\infty}_{[s,T]}})\big)\sum_{i=1}^{n}\exp(\gamma|Y_{s}^{i}|) {\rm d}s.
\end{array}
\end{equation} 
Noting that 
$$
\frac{\exp(\gamma|x|)-2}{2\gamma^2} \leq u(x) \leq \frac{\exp(\gamma|x|)}{\gamma^2}.
$$ 
We have
\begin{equation}\label{eq:3.8}
\begin{array}{lll}
\Dis \sum_{i=1}^{n}\frac{\exp(\gamma|Y_{t}^{i}|)-2}{2\gamma^2} & \leq & \Dis  \frac{n\exp(\gamma\|\xi\|_{\infty})}{\gamma^2}-\int_{t}^{T} \sum_{i=1}^{n}u'(Y_{s}^{i})Z_{s}^{i}{\rm d}W_s \vspace{0.2cm}\\
& & \Dis +\frac{1}{\gamma}\int_{t}^{T}\big(\alpha_{s}+\beta_{s}\|Y\|_{\s^{\infty}_{[s,T]}}+k_s\exp(\frac{r\gamma}{2-r}\|Y\|_{\s^{\infty}_{[s,T]}})\big)\sum_{i=1}^{n}\exp(\gamma|Y_{s}^{i}|) {\rm d}s.
\end{array}
\end{equation} 
Taking conditional expectation with respect to $\F_{\tau}$ for $\tau \in [t_0,t]$, we show that
\begin{equation}\label{eq:3.9}
\begin{array}{lll}
&& \Dis \E\big[\sum_{i=1}^{n}\exp(\gamma|Y_{t}^{i}|)|\F_{\tau}\big] \vspace{0.2cm}\\
& \leq & \Dis 2n(\exp(\gamma\|\xi\|_{\infty})+1) \vspace{0.2cm}\\
&& \Dis +\int_{t}^{T}2\gamma\big(\alpha_{s}+\beta_{s}\|Y\|_{\s^{\infty}_{[s,T]}}+k_s\exp(\frac{r\gamma}{2-r}\|Y\|_{\s^{\infty}_{[s,T]}})\big)\E\big[\sum_{i=1}^{n}\exp(\gamma|Y_{s}^{i}|)|\F_{\tau}\big] {\rm d}s.
\end{array}
\end{equation}  
Using Gronwall's inequality, we get 
\begin{equation}\label{eq:3.10}
\begin{array}{lll}
&&\Dis \E\big[\sum_{i=1}^{n}\exp(\gamma|Y_{t}^{i}|)|\F_{\tau}\big] \vspace{0.2cm}\\
& \leq & \Dis  2n(\exp(\gamma\|\xi\|_{\infty})+1)\cdot\exp\Big(\int_{t}^{T}2\gamma\big(\alpha_{s}+\beta_{s}\|Y\|_{\s^{\infty}_{[s,T]}}+k_s\exp(\frac{r\gamma}{2-r}\|Y\|_{\s^{\infty}_{[s,T]}})\big){\rm d}s\Big).
\end{array}
\end{equation} 
Setting $\tau=t$ and noting $r \leq 1+\delta$, we have
$$
\sum_{i=1}^{n}\exp(\gamma|Y_{t}^{i}|)  \leq   2n\big(\exp(\gamma\|\xi\|_{\infty})+1\big)\exp\Big(\int_{t}^{T}2\gamma\big(\alpha_{s}+\beta_{s}\|Y\|_{\s^{\infty}_{[s,T]}}+k_s\exp(\frac{r\gamma}{1-\delta}\|Y\|_{\s^{\infty}_{[s,T]}})\big){\rm d}s\Big).
$$
Using Jensen's inequality, we obtain that 
$$
\sum_{i=1}^{n}\exp(\gamma|Y_{t}^{i}|) \geq n\exp(\frac{\sum_{i=1}^{n}\gamma|Y_{t}^{i}|}{n}) \geq n\exp(\frac{\gamma|Y_{t}|}{n}).
$$
Combining the preceding inequalities and the assumption \ref{A:H3}, we have
\begin{equation}\label{eq:3.11}
\begin{array}{lll}
|Y_{t}| & \leq & \Dis \frac{n}{\gamma}\log(2\exp(\gamma\|\xi\|_{\infty})+2)+ \int_{t}^{T}2n\big(\alpha_{s}+\beta_{s}\|Y\|_{\s^{\infty}_{[s,T]}}+k_s\exp(\frac{r\gamma}{1-\delta}\|Y\|_{\s^{\infty}_{[s,T]}})\big){\rm d}s \vspace{0.2cm}\\
& \leq & \Dis \frac{n}{\gamma}\log(2\exp(\gamma C_0)+2)+ 2nC_1+2n\int_{t_0}^{T} k_s\exp(\frac{r\gamma}{1-\delta}\|Y\|_{\s^{\infty}_{[s,T]}})\big){\rm d}s \\
&& \Dis + \int_{t}^{T}2n\beta_{s}\|Y\|_{\s^{\infty}_{[s,T]}}{\rm d}s.
\end{array}
\end{equation} 
Let 
$$
K_0 :=\frac{n}{\gamma}\log(2\exp(\gamma C_0)+2)+ 2nC_1+2n\int_{t_0}^{T} k_s\exp(\frac{r\gamma}{1-\delta}\|Y\|_{\s^{\infty}_{[s,T]}})\big){\rm d}s.
$$
We have
\begin{equation}\label{eq:3.12}
\|Y\|_{\s^{\infty}_{[t,T]}} \leq K_0+\int_{t}^{T}2n\beta_{s}\|Y\|_{\s^{\infty}_{[s,T]}}{\rm d}s, \ \ \forall t \in [t_0,T].
\end{equation}
Using Gronwall's inequality and the assumption \ref{A:H3}, we have
\begin{equation}\label{eq:3.13}
\begin{array}{lll}
\|Y\|_{\s^{\infty}_{[t_0,T]}} &\leq& \Dis K_{0}\exp(\int_{t_0}^{T}2n\beta_{s}{\rm d}s)\vspace{0.2cm}\\ 
&\leq& \Dis \exp(2nC_2)\bigg(\frac{n}{\gamma}\log(2\exp(\gamma C_0)+2)+ 2nC_1+2n\int_{t_0}^{T} k_s\exp(\frac{r\gamma}{1-\delta}\|Y\|_{\s^{\infty}_{[s,T]}})\big){\rm d}s\bigg).
\end{array}
\end{equation}
Let 
$$
C_4=\exp(2nC_2)\Big(\frac{n}{\gamma}\log\big(2\exp(\gamma C_0)+2\big)+ 2nC_1\Big), \ \ C_5=2nC_3\exp(2nC_2)\Big((\frac{2n}{\gamma})^{\frac{1+\delta}{1-\delta}}+1\Big).
$$
From the definition of $k_s$ and the assumption \ref{A:H3}, we get \eqref{eq:3.1}. The proof is complete.\vspace{0.2cm}

From Lemma~\ref{lem:3.1}, we get the following proposition.

\begin{prop}\label{pro:3.2}
There exists a constant $r_0>0$ (depending on the vector of parameters $(n,\gamma,\delta,C_0,C_1,C_2,C_3)$) such that if \ref{A:H1}-\ref{A:H3} holds for $r \in (0,r_0)$, and $(Y,Z) \in \s^{\infty}_{[t_0,T]}(\R^n) \times {\hcal^2}_{[t_0,T]}(\R^{n\times d})$ is a solution of BSDE \eqref{eq:1.1} on $[t_0,T]$, then 
\begin{equation}\label{eq:3.14}
\|Y\|_{\s^{\infty}_{[t_0,T]}} \leq C_4+2C_5,
\end{equation}
where $C_4$ and $C_5$ are given by \eqref{eq:3.1}.
\end{prop}

{\bf Proof}\ \ Define
$$
F(x)=C_4+C_5\exp(\frac{r\gamma x}{1-\delta})-x , \ \ x \geq 0.
$$
Then we have
$$
F^{'}(x)=\frac{r\gamma C_5}{1-\delta}\exp(\frac{r\gamma x}{1-\delta})-1, \ \ F^{''}(x)=\frac{r^{2}\gamma^{2} C_5}{(1-\delta)^{2}}\exp(\frac{r\gamma x}{1-\delta}) > 0.
$$
Let 
$$
r_0=\frac{(1-\delta)\log 2}{\gamma(C_4+2C_5)}.
$$ 
For a given $r \in (0,r_0)$, let 
$$
x_0=\frac{1-\delta}{r\gamma}\log \frac{1-\delta}{C_{5}r\gamma}.
$$ 
Then we have 
$$
r < r_0 < \frac{1-\delta}{\gamma C_5}, \ \ x_0 > 0, \ \ F^{'}(x_0) = 0.
$$
Hence $F$ is decreasing on $[0,x_0]$ and increasing on $[x_0,+\infty)$, and$\vspace{0.2cm}$
\begin{equation*}
\begin{array}{lll}
\Dis F(C_4+2C_5) & = & \Dis C_4+C_5\exp(\frac{r\gamma(C_4+2C_5)}{1-\delta})-C_4-2C_5 \vspace{0.2cm}\\
& = & \Dis C_5(\exp(\frac{r\gamma(C_4+2C_5)}{1-\delta})-2) \vspace{0.2cm}\\
& < & \Dis C_5(\exp(\frac{r_0\gamma(C_4+2C_5)}{1-\delta})-2)=0.
\end{array}
\end{equation*}
Then $F(x)=0$ has two zeros $x_1,x_2$ and they satisfy
\begin{equation}\label{eq:3.15}
C_0 < C_4 < x_1 < C_4+2C_5 < x_2, \ \ \{x:F(x)\geq 0\}=[0,x_1]\cup [x_2,+\infty).
\end{equation}
From Lemma~\ref{lem:3.1}, we obtain that 
$$
F(\|Y\|_{\s^{\infty}_{[t,T]}}) \geq 0,\ \forall t \in [t_0,T].
$$ 
Hence 
$$
\|Y\|_{\s^{\infty}_{[t,T]}} \in [0,x_1]\cup [x_2,+\infty),\ \forall t \in [t_0,T].
$$ 
Define
$$
\hat{t}=\inf\{t\in [t_0,T]:\|Y\|_{\s^{\infty}_{[t,T]}}\leq x_1\}.
$$
Notice that 
$$
\|Y\|_{\s^{\infty}_{[T,T]}}=\|\xi\|_{\infty}\leq C_0 < x_1.
$$
$\hat{t}$ is well defined. $\|Y\|_{\s^{\infty}_{[t,T]}}$ is decreasing and right-continuous about $t$, so we have $\|Y\|_{\s^{\infty}_{[\hat{t},T]}} \leq x_1$. If $\hat{t} > t_0$, then
$$
\|Y\|_{\s^{\infty}_{[t,T]}} > x_1, \ \ \forall t \in [t_0,\hat{t}).
$$
Therefore,
$$
\|Y\|_{\s^{\infty}_{[t,T]}} \geq x_2, \ \ \forall t \in [t_0,\hat{t}).
$$
From (\ref{eq:3.13}), we deduce that 
\begin{equation}\label{eq:3.16}
\begin{array}{lll}
\Dis x_2 & \leq & \Dis \limsup_{t\to \hat{t}^{-}}\|Y\|_{\s^{\infty}_{[t,T]}} \vspace{0.2cm}\\
& \leq & \Dis \limsup_{t\to \hat{t}^{-}}\Big[C_4+2n\exp(2nC_2) \int_{t}^{T}k_s\exp(\frac{r\gamma}{1-\delta}\|Y\|_{\s^{\infty}_{[s,T]}}){\rm d}s\Big] \vspace{0.2cm}\\
& = & \Dis C_4+ 2n\exp(2nC_2)\int_{\hat{t}}^{T}k_s\exp(\frac{r\gamma}{1-\delta}\|Y\|_{\s^{\infty}_{[s,T]}}){\rm d}s \vspace{0.2cm}\\
& \leq & \Dis C_4+2n\exp(2nC_2)((\frac{2n}{\gamma})^{\frac{1+\delta}{1-\delta}}+1)C_3\exp(\frac{r\gamma}{1-\delta}\|Y\|_{\s^{\infty}_{[\hat{t},T]}}) \vspace{0.2cm}\\
& \leq & \Dis C_4+C_5\exp(\frac{r\gamma x_1}{1-\delta}) = F(x_1)+x_1=x_1.
\end{array}
\end{equation} 
This is a contradiction. Hence $\hat{t}=t_0$, and
$$
\|Y\|_{\s^{\infty}_{[t_0,T]}} \leq x_1 < C_4+2C_5.
$$
The proof is complete.\vspace{0.2cm}

{\bf Proof of Theorem~\ref{thm:2.1}}\ \ For the number $r_0$ given in Proposition~\ref{pro:3.2} and a given $r \in (0,r_0)$, define
$$
\lambda := C_4+2C_5,
$$
where $C_4$ and $C_5$ are the same as in Lemma~\ref{lem:3.1}. From \eqref{eq:3.15}, we have 
$$
\|\xi\|_{\infty} \leq C_0 \leq C_4 \leq \lambda.
$$
From \cite[Theorem 2.2, p. 1072]{HuTang2016SPA}, there exists $t_{\lambda} > 0$ which depends on constants $(n,C,\gamma,\delta,\lambda)$, such that BSDE \eqref{eq:1.1} has a local solution $(Y,Z) \in \s^{\infty}(\R^n) \times {\rm BMO}(\R^{n\times d}) $ on $[T-t_{\lambda},T]$. From Proposition~\ref{pro:3.2}, we obtain that 
$$
\|Y_{T-t_{\lambda}}\|_{\infty} \leq \|Y\|_{\s^{\infty}_{[T-t_{\lambda},T]}} \leq \lambda.
$$
Taking $T-t_{\lambda}$ as the terminal time and $Y_{T-t_{\lambda}}$ as the terminal value, BSDE \eqref{eq:1.1} has a local solution $(Y,Z) \in \s^{\infty}(\R^n) \times {\rm BMO}(\R^{n\times d}) $ on $[T-2t_{\lambda},T-t_{\lambda}]$. Stitching the solutions we have a solution $(Y,Z) \in \s^{\infty}(\R^n) \times {\hcal^2}(\R^{n\times d})$ on $[T-2t_{\lambda},T]$ and $\|Y_{T-2t_{\lambda}}\|_{\infty} \leq \lambda$. Repeating the preceding process, we can extend the pair $(Y,Z)$ to the whole interval $[0,T]$ within finite steps such that $Y$ is uniformly bounded by $\lambda$ and $Z \in {\hcal^2}(\R^{n\times d})$. We now show that  $Z \in {\rm BMO}(\R^{n\times d})$.
Identical to the proof of inequality \eqref{eq:3.5} and \eqref{eq:3.6}, we have
\begin{equation}\label{eq:3.17}
\Dis \eta_{s}|Z_{s}|^{r} \leq \frac{\gamma}{4n}\exp(-\gamma\|Y\|_{\s^{\infty}_{[s,T]}})|Z_{s}|^{2} + (\eta_{s}^{\frac{2}{1-\delta}}+\eta_{s})((\frac{4n}{\gamma})^{\frac{1+\delta}{1-\delta}}+1)\exp(\frac{r\gamma}{2-r}\|Y\|_{\s^{\infty}_{[s,T]}}),
\end{equation}
and
\begin{equation}\label{eq:3.18}
\begin{array}{lll}
\Dis u(Y_{t}^{i}) & \leq & \Dis  u(\xi^{i})-\int_{t}^{T}u'(Y_{s}^{i})Z_{s}^{i}{\rm d}W_s+\int_{t}^{T}\big[-\frac{1}{2}|Z_{s}^{i}|^{2}+\frac{1}{4n}|Z_{s}|^{2}\big] {\rm d}s \vspace{0.2cm}\\
& & \Dis +\int_{t}^{T}\frac{\exp(\gamma|Y_{s}^{i}|)}{\gamma}\big(\alpha_{s}+\beta_{s}\|Y\|_{\s^{\infty}_{[s,T]}}+\hat{k}_s\exp(\frac{r\gamma}{2-r}\|Y\|_{\s^{\infty}_{[s,T]}})\big) {\rm d}s,
\end{array}
\end{equation}
where 
$$
\hat{k}_s=(\eta_{s}^{\frac{2}{1-\delta}}+\eta_{s})((\frac{4n}{\gamma})^{\frac{1+\delta}{1-\delta}}+1).
$$
Summing $i$ from $1$ to $n$ and taking conditional expectation with respect to $\F_{t}$, we have
\begin{equation*}
\Dis \frac{1}{4}\E\big[\int_{t}^{T}|Z_{s}|^{2}{\rm d}s|\F_{t}\big] \leq \frac{n\exp(\gamma C_0)}{\gamma^2} + \frac{n\exp(\gamma \lambda)}{\gamma}\bigg(C_{1}+C_{2}\lambda +C_{3}\big((\frac{4n}{\gamma})^{\frac{1+\delta}{1-\delta}}+1\big)\exp\big(\frac{r\gamma\lambda}{2-r}\big)\bigg).
\end{equation*}
Hence $Z \in {\rm BMO}(\R^{n\times d})$. Finally, the uniqueness on the given interval $[0,T]$ is a consequence of \cite[Theorem 2.2, p. 1072]{HuTang2016SPA} via a pasting technique.
\begin{rem}\label{rmk:3.1}
Assumptions \ref{A:H1}-\ref{A:H3} of Theorem~\ref{thm:2.1} are different from those of \cite[Theorem 2.3, p. 1072]{HuTang2016SPA} and \cite[Theorem 2.4]{FanHuTang2020ArXiv}. We allow the generator to have a small growth of the off-diagonal elements. They are different from those of \cite[Theorem 2.5]{FanHuTang2020ArXiv} in that the generator is not required to be strictly quadratic. For example, the following generator $f$ satisfies Theorem~\ref{thm:2.1} rather than the others when $r$ is sufficiently small: $$f^i(t,y,z)=|z^i|^2\sin\big(\log(|z^i|+1)\big)+|y|+\sin(|z|^{1+\delta})+|z|^{r}{\bf 1}_{\{|z|>1\}}+|z|{\bf 1}_{\{|z|\leq 1 \}},\ \ i=1,\cdots, n.$$
\end{rem}
\begin{rem}\label{rmk:3.2}
When $C_3$ is sufficiently small such that $C_5<\exp(-\frac{\gamma(1+\delta)(1+C_4)}{1-\delta})$, taking $r_0=1+\delta$, then for $r \in (0,r_0]$, we have
$$
F(C_4+1)=C_5\exp(\frac{\gamma(1+\delta)(1+C_4)}{1-\delta})-1 <0.
$$
In a similar way we have Theorem~\ref{thm:2.1}. In particular, when $C_3=0$ we have $C_5=0$, then we have $\|Y\|_{\s^{\infty}_{[t_0,T]}} \leq C_4$ by \eqref{eq:3.1} and thus Theorem~\ref{thm:2.1} holds, which is the case of \cite[Theorem 2.4]{FanHuTang2020ArXiv}.
\end{rem}
\begin{rem}\label{rmk:3.3}
From \cite[Theorem 2.1]{FanHuTang2020ArXiv}, \ref{A:H1} and \ref{A:H2} can be replaced with the following in Theorem~\ref{thm:2.1}:
\begin{enumerate}
\renewcommand{\theenumi}{(H\arabic{enumi}')}
\renewcommand{\labelenumi}{\theenumi}
\item\label{A:H1'} There exist a deterministic scalar-valued positive function $(\alpha_{t})_{t\in\T}$, a deterministic nondecreasing continuous function $\phi(\cdot):[0,+\infty)\To [0,+\infty)$ with $\phi(0)=0$ and several real constants $\gamma>0$, $C \geq 0$, $\delta\in [0,1)$ such that for $i=1,\cdots,n$ and each $(y,\bar y, z, \bar z)\in \R^n\times\R^n\times\R^{n\times d}\times\R^{n\times d}$, $f^i$ satisfies the following inequalities:
    $$
    \begin{array}{ll}
    &\Dis |f^i(t,y,z)|\leq \alpha_t+\phi(|y|)+\frac{\gamma}{2} |z^i|^2+C\sum_{j\neq i} |z^j|^{1+\delta}; \vspace{0.1cm}\\
    &\Dis |f^i(t,y,z)-f^i(t,\bar y,\bar z)| \vspace{0.1cm}\\
    &\Dis \leq \phi(|y|\vee|\bar y|) \Dis \left[ \left(1+|z|+|\bar z|\right)\left(|y-\bar y|+|z^i-\bar z^i|\right)+\left(1+|z|^\delta+|\bar z|^\delta\right)\sum_{j\neq i} |z^j-\bar z^j|\right].
    \end{array}
    $$
\item\label{A:H2'} There exist a two-dimensional non-negative deterministic vector function $(\beta_{t},\eta_{t})_{t\in\T}$ and a positive constant $r \in (0,1+\delta]$ such that for $i=1,\cdots,n$ and $(y,z)\in \R^n\times\R^{n\times d}$, the function $f^i$ satisfies:
    $$
    {\rm sgn}(y^i)f^i(t,y,z)\leq \alpha_{t}+\beta_{t}|y|+\eta_{t}|z|^{r}+\frac{\gamma}{2} |z^i|^2.
    $$
\end{enumerate}
\end{rem}

\Section{Diagonally quadratic and triangular BSDEs}
\vspace{0mm}

To prove Theorem~\ref{thm:2.2}, we need the following lemma.
\begin{lem}\label{lem:4.1}
We consider the following one-dimensional BSDE:
\begin{equation}\label{eq:4.1}
\Dis Y_{t}=\eta+\int_{t}^{T}l(s,Y_{s},Z_{s}){\rm d}s-\int_{t}^{T}Z_{s}{\rm d}W_s,\ 0 \leq t \leq T.
\end{equation}
The terminal value $\eta$ and the generator $l$ satisfy the following assumptions:
\begin{enumerate}
\renewcommand{\theenumi}{(B\arabic{enumi})}
\renewcommand{\labelenumi}{\theenumi}
\item\label{B:B1} $\forall ~K>0$, the function $l:\Omega\times\T\times\R\times\R^{1\times d}\To \R$ satisfies:
    $$
    \sup_{t\in\T}\left\|\E\Big[\exp\Big(K\int_{t}^{T}\big|l(s,0,0)\big|{\rm d}s\Big)\Big|\F_t\Big]\right\|_{\infty} < +\infty.
    $$

\item\label{B:B2} There exist a non-negative constant $\beta$ and a positive constant $C$ such that for each $(y,\bar y, z, \bar z)\in \R\times\R\times\R^{1\times d}\times\R^{1\times d}$, the function $l$ satisfies:
    $$
    |l(t,y,z)-l(t,\bar y,\bar z)| \leq \beta|y-\bar y|+C(1+|z|+|\bar z|)|z-\bar z|.
    $$

\item\label{B:B3} There exists a non-negative constant $C_1$ such that $\eta$ satisfies:
$$\|\eta\|_{\infty}\leq C_1.$$

\end{enumerate}
Then BSDE \eqref{eq:4.1} has a unique solution $(Y,Z)\in \s^{\infty}(\R) \times {\rm BMO}(\R^{1\times d})$ on $[0,T]$.
\end{lem}
{\bf Proof}\ \ When $\beta=0$, $l$ is independent of $y$. From \cite[Lemma 2.1]{HuTang2016SPA}, we know the result holds. When $\beta > 0$, for $y\in \s^{\infty}(\R)$, we define a map $\varphi(y)=Y$, where $Y$ is given by:
\begin{equation*}
\Dis Y_{t}=\eta+\int_{t}^{T}l(s,y_{s},Z_{s}){\rm d}s-\int_{t}^{T}Z_{s}{\rm d}W_s,\ 0 \leq t \leq T.
\end{equation*}
From the preceding result, we know $\varphi$ is well-defined and maps $\s^{\infty}(\R)$ to itself. For $y,\bar y \in \s^{\infty}(\R)$, let $Y=\varphi(y),\bar Y=\varphi(\bar y)$. Denote $\delta Y=Y-\bar Y, \delta y=y-\bar y,\delta Z=Z-\bar Z$. We have
\begin{equation*}
\begin{array}{lll}
\Dis \delta Y_{t}&=&\int_{t}^{T}[l(s,y_{s},Z_{s})-l(s,\bar y_{s},\bar Z_{s})]{\rm d}s-\int_{t}^{T}\delta Z_{s}{\rm d}W_s\vspace{0.2cm}\\
&=& \Dis \int_{t}^{T}(\beta_{s}\delta y_{s}+\delta Z_{s}\alpha_{s}){\rm d}s-\int_{t}^{T}\delta Z_{s}{\rm d}W_s.
\end{array}
\end{equation*}
Here $|\beta_{s}| \leq \beta, ~|\alpha_{s}| \leq C(1+|Z_{s}|+|\bar Z_{s}|)$, therefore $\alpha\cdot W$ is a BMO martingale. Define
\begin{equation*}
\widetilde{W}_t:=W_t-\int_0^t\alpha_s\, {\rm d}s, \ \ \ t\in [0,T]; \ \ \ {\rm d}\widetilde {\mathbb{P}}:=\mathscr{E}(\alpha\cdot W)_0^T {\rm d}\mathbb{P}.
\end{equation*}
Then, $\widetilde {\mathbb{P}}$ is a new probability equivalent to $\mathbb{P}$, and $\widetilde W$ is a Brownian motion with respect to $\widetilde {\mathbb{P}}$. We have
\begin{equation*}
\Dis \delta Y_{t}=\int_{t}^{T}\beta_{s}\delta y_{s}{\rm d}s-\int_{t}^{T}\delta Z_{s}{\rm d}\widetilde{W}_s.
\end{equation*}
Taking the conditional expectation with respect to $\widetilde {\mathbb{P}}$, we have
\begin{equation*}
\Dis \|\delta Y\|_{{\s}^{\infty}} \leq \int_{t}^{T}\beta\|\delta y\|_{{\s}^{\infty}}{\rm d}s \leq \beta T \|\delta y\|_{{\s}^{\infty}}.
\end{equation*}
When $T \leq \frac{1}{2\beta}$, $\varphi$ is a contraction map and and the statement follows from the Banach fixed point theorem. For general $T$, we can repeat the preceding process and get the result within finite steps. The proof is complete.\vspace{0.2cm}

{\bf Proof of Theorem~\ref{thm:2.2}}\ \ We will solve BSDE \eqref{eq:2.2} in order. For the first equation, notice that $\|\xi^{1}\|_{\infty} \leq K_3$, $|k^{1}(s,0,0)| \leq K_1$, from \cite{Kobylanski2000AP} we know it has a unique solution $(Y^{1},Z^{1})\in \s^{\infty}(\R) \times {\rm BMO}(\R^{1\times d})$. Suppose that we already solve the first $(i-1)$ equations with $(Y^{j},Z^{j})\in \s^{\infty}(\R) \times {\rm BMO}(\R^{1\times d}),j=1,\cdots,i-1$. For the $i$-th equation, we have $\forall K > 0$,
\begin{equation}\label{eq:4.2}
\begin{array}{lll}
\Dis & & \Dis \E\Big[\exp\Big(K\int_{t}^{T}\big|k^{i}(s,Y_{s}(0;i),Z_{s}(0;i))\big|{\rm d}s\Big)\Big|\F_t\Big]\vspace{0.2cm}\\
&\leq& \Dis \exp\big(KTK_1\big)\cdot\E\Big[\exp\Big(KK_1\int_{t}^{T}(\sum_{j=1}^{i-1}|Y_{s}^j|+\sum_{j=1}^{i-1}|Z_{s}^j|^{1+\alpha}){\rm d}s\Big)\Big|\F_{t}\Big] \vspace{0.2cm}\\
&\leq& \Dis \exp\Big(KTK_1+KTK_1\sum_{j=1}^{i-1}\|Y^j\|_{\s^{\infty}}\Big)\cdot\E\Big[\exp\Big(KK_1\int_{t}^{T}\sum_{j=1}^{i-1}|Z_{s}^j|^{1+\alpha}{\rm d}s\Big)\Big|\F_{t}\Big]. 
\end{array}
\end{equation}
By H\"{o}lder's inequality and Young's inequality we get 
\begin{equation}\label{eq:4.3}
\E\Big[\exp\Big(KK_1\int_{t}^{T}\sum_{j=1}^{i-1}|Z_{s}^j|^{1+\alpha}{\rm d}s\Big)\Big|\F_{t}\Big] \leq \bigg(\prod_{j=1}^{i-1}\E\Big[\exp\Big(KK_1(i-1)\int_{t}^{T}|Z_{s}^j|^{1+\alpha}{\rm d}s\Big)\Big|\F_{t}\Big]\bigg)^{\frac{1}{i-1}},
\end{equation}
and
$$
L|Z_{s}^j|^{1+\alpha}=\varepsilon^{\frac{1+\alpha}{2}}|Z_{s}^j|^{1+\alpha}\cdot L\varepsilon^{-\frac{1+\alpha}{2}} \leq \frac{1+\alpha}{2}\varepsilon|Z_{s}^j|^{2}+\frac{1-\alpha}{2}L^{\frac{2}{1-\alpha}}\varepsilon^{-\frac{1+\alpha}{1-\alpha}},\ \ \forall L >0.
$$
Let $\varepsilon$ be sufficiently small such that $\frac{1+\alpha}{2}\varepsilon\|Z^j\|_{\rm BMO}^{2} \leq \frac{1}{2}$. From John-Nirenberg inequality~(\cite[Theorem 2.2]{Kazamaki1994book}), we have for $1 \leq j \leq i-1$,
\begin{equation}\label{eq:4.4}
\begin{array}{lll}
\Dis \E\Big[\exp\Big(L\int_{t}^{T}|Z_{s}^j|^{1+\alpha}{\rm d}s\Big)\Big|\F_{t}\Big] &\leq& \Dis \frac{1}{1-\frac{1+\alpha}{2}\varepsilon\|Z^j\|_{\rm BMO}^{2}}\exp\Big(\frac{1-\alpha}{2}L^{\frac{2}{1-\alpha}}\varepsilon^{-\frac{1+\alpha}{1-\alpha}}T\Big) \vspace{0.4cm}\\
&\leq& \Dis 2\exp\Big(\frac{1-\alpha}{2}L^{\frac{2}{1-\alpha}}\varepsilon^{-\frac{1+\alpha}{1-\alpha}}T\Big), \ \forall ~t \in\T.
\end{array}
\end{equation}
Combining (\ref{eq:4.2}),(\ref{eq:4.3}) and (\ref{eq:4.4}), we obtain that $\forall ~K>0$,
\begin{equation}\label{eq:4.5}
\sup_{t\in\T}\left\|\E\Big[\exp\Big(K\int_{t}^{T}\big|k^{i}(s,Y_{s}(0;i),Z_{s}(0;i))\big|{\rm d}s\Big)\Big|\F_t\Big]\right\|_{\infty} < +\infty.
\end{equation}
Therefore, we can apply Lemma~\ref{lem:4.1} to see the $i$-th equation admits a unique solution $(Y^{i},Z^{i})\in \s^{\infty}(\R) \times {\rm BMO}(\R^{1\times d})$ on $[0,T]$. The proof is complete.

\Section{Conclusion Remark}
\vspace{0mm}

We study the well-posedness of the multi-dimensional BSDE~(\ref{eq:1.1}) with a diagonally quadratic generator. When the generator has a small growth of the off-diagonal elements, we build a new priori estimate and get the existence and uniqueness of the global solution, which generalizes the results in Hu and Tang~\cite{HuTang2016SPA} and Fan et al.~\cite{FanHuTang2020ArXiv}. Besides, when the generator is diagonally quadratic and triangular, we get the global solvability of the multi-dimensional BSDE~(\ref{eq:2.2}) without the small growth condition. Finally, when the generator is non-triangular and has a general sub-quadratic growth of the off-diagonal elements, the existence and uniqueness of the global solutions are interesting and challenging, which remains to be studied in the future.

\end{document}